\documentclass[11pt]{article}
\usepackage{amsmath,amssymb}
\usepackage{graphicx, cite}

\newtheorem{theorem}{Theorem}

\newcommand{\<}{\langle}
\renewcommand{\>}{\rangle}
\newcommand{\R}{\mathbb{R}}

\newcommand{\ncol}{m}
\newcommand{\nrow}{n}

\begin{document}

\title{Stable Signal Recovery from \\
Incomplete and Inaccurate Measurements}

\author{Emmanuel Candes$^{\dagger}$, Justin Romberg$^{\dagger}$,
  and Terence Tao$^{\sharp}$\\
  \vspace{-.2cm}\\
  $\dagger$ Applied and Computational Mathematics, 
Caltech, Pasadena, CA 91125\\
  \vspace{-.3cm}\\
  $\sharp$ Department of Mathematics, University of California, Los
  Angeles, CA 90095 }

\date{February, 2005; Revised June 2005} 

\maketitle
 
\begin{abstract}
  Suppose we wish to recover a vector $x_0 \in \R^\ncol$ (e.g. a
  digital signal or image) from incomplete and contaminated
  observations $y = A x_0 + e$; $A$ is a $\nrow$ by $\ncol$ matrix
  with far fewer rows than columns ($\nrow \ll \ncol$) and $e$ is an
  error term.  Is it possible to recover $x_0$ accurately based on the
  data $y$?
  
  To recover $x_0$, we consider the solution $x^\sharp$ to the
  $\ell_1$-regularization problem
\[
\min~\|x\|_{\ell_1} \quad\text{subject to} \quad 
\|Ax-y\|_{\ell_2} \leq \epsilon,
\]
where $\epsilon$ is the size of the error term $e$. We show that if
$A$ obeys a uniform uncertainty principle (with unit-normed columns)
and if the vector $x_0$ is sufficiently sparse, then the solution is
within the noise level
\[
\|x^\sharp - x_0\|_{\ell_2} \le C \cdot \epsilon. 
\]
As a first example, suppose that $A$ is a Gaussian random matrix, then
stable recovery occurs for almost all such $A$'s provided that the
number of nonzeros of $x_0$ is of about the same order as the number
of observations.  As a second instance, suppose one observes few
Fourier samples of $x_0$, then stable recovery occurs for almost any
set of $\nrow$ coefficients provided that the number of nonzeros is of
the order of $\nrow/[\log \ncol]^6$.

In the case where the error term vanishes, the recovery is of course
exact, and this work actually provides novel insights on the exact
recovery phenomenon discussed in earlier papers. The methodology also
explains why one can also very nearly recover approximately sparse
signals.
\end{abstract}


{\bf Acknowledgments.} E.~C. is partially supported by a National
Science Foundation grant DMS 01-40698 (FRG) and by an Alfred P.  Sloan
Fellowship.  J.~R.~is supported by National Science Foundation grants
DMS 01-40698 and ITR ACI-0204932.  T.~T.~is supported in part by
grants from the Packard Foundation.  

\pagebreak

%----------------------------------------------------------------------
\section{Introduction}
\label{sec:intro}

\subsection{Exact recovery of sparse signals}

Recent papers
\cite{candes04ro,candes04qu,candes04ne,candes04de,donoho04fo} have
developed a series of powerful results about the exact recovery of a
finite signal $x_0 \in \R^\ncol$ from a very limited number of
observations. As a representative result from this literature,
consider the problem of recovering an unknown {\em sparse} signal
$x_0(t)\in\R^\ncol$; that is, a signal $x_0$ whose support $T_0 = \{t
: x_0(t)\not=0\}$ is assumed to have small cardinality.  All we know
about $x_0$ are $\nrow$ linear measurements of the form
\[
y_k = \<x_0,a_k\> \quad k = 1,\ldots,\nrow  \quad\text{or}\quad y = Ax_0, 
\]   
where the $a_k\in\R^\ncol$ are known test signals. Of special interest
is the vastly underdetermined case, $\nrow \ll \ncol$, where there are
many more unknowns than observations. At first glance, this may seem
impossible.  However, it turns out that one can actually recover $x_0$
exactly by solving the convex program\footnote{$(P_1)$ can even be
  recast as a linear program \cite{chen99at}.}
\begin{equation}
\label{eq:P1}
(P_1)\quad\quad \min~\|x\|_{\ell_1} \quad\text{subject to}\quad Ax=y,
\end{equation}
provided that the matrix $A \in \R^{\nrow \times \ncol}$ obeys a
{\em uniform uncertainty principle}.

The uniform uncertainty principle, introduced in \cite{candes04ne} and
refined in \cite{candes04de}, essentially states that the $\nrow
\times \ncol$ measurement matrix $A$ obeys a ``restricted isometry
hypothesis.''  To introduce this notion, let $A_T,~T\subset \{1,
\ldots, \ncol\}$ be the $\nrow \times |T|$ submatrix obtained by
extracting the columns of $A$ corresponding to the indices in $T$.
Then \cite{candes04de} defines the $S$-restricted isometry constant
$\delta_S$ of $A$ which is the smallest quantity such that
\begin{equation}
\label{eq:riso}
(1-\delta_S)\,\|c\|^2_{\ell_2} \leq \|A_T c\|^2_{\ell_2}  \leq 
(1+\delta_S)\, \|c\|^2_{\ell_2}
\end{equation}
for all subsets $T$ with $|T|\leq S$ and coefficient sequences
$(c_j)_{j\in T}$. This property essentially requires that every set of
columns with cardinality less than $S$ approximately behaves like an
orthonormal system. 
It was shown (also in \cite{candes04de}) that if $S$ verifies 
\begin{equation}
\label{eq:d3s}
\delta_S + \delta_{2S} + \delta_{3S} < 1, 
\end{equation}
then solving $(P_1)$ recovers {\em any} sparse signal $x_0$ with
support size obeying $|T_0|\leq S$.

%----------------------------------------------------------------------
\subsection{Stable recovery from imperfect measurements}

This paper develops results for the ``imperfect'' (and far more
realistic) scenarios where the measurements are noisy and the signal
is not exactly sparse.  Everyone would agree that in most practical
situations, we cannot assume that $A x_0$ is known with arbitrary
precision. More appropriately, we will assume instead that one is
given ``noisy'' data $y = A x_0 + e$, where $e$ is some unknown
perturbation bounded by a known amount $\|e\|_{\ell_2} \leq
\epsilon$. To be broadly applicable, our recovery procedure must be
{\em stable}: small changes in the observations should result in small
changes in the recovery.  This wish, however, may be quite hopeless.
How can we possibly hope to recover our signal when not only the
available information is severely incomplete but in addition, the few
available observations are also inaccurate?

Consider nevertheless (as in \cite{donoho04st} for example) the convex
program searching, among all signals consistent with the data $y$, for
that with minimum $\ell_1$-norm
\begin{equation}
\label{eq:P2}
(P_2) \quad\quad 
\min~\|x\|_{\ell_1} \quad\text{subject to}\quad \|Ax-y\|_{\ell_2} 
\leq \epsilon.
\end{equation}
The first result of this paper shows that contrary to the belief
expressed above, the solution to $(P_2)$ recovers an unknown sparse
object with an error at most proportional to the noise level.  Our
condition for stable recovery again involves the restricted isometry
constants.
\begin{theorem}
\label{th:stable}
Let $S$ be such that $\delta_{3S} + 3 \delta_{4S} < 2$.  Then for any
signal $x_0$ supported on $T_0$ with $|T_0|\leq S$ and any
perturbation $e$ with $\|e\|_{\ell_2} \leq \epsilon$, the solution
$x^\sharp$ to $(P_2)$ obeys
\begin{equation}
\label{eq:lqerror}
\|x^\sharp - x_0 \|_{\ell_2} \leq C_S\cdot\epsilon,
\end{equation}
where the constant $C_S$ may only depend on $\delta_{4S}$.  For
reasonable values of $\delta_{4S}$, $C_S$ is well behaved; e.g. $C_S
\approx 8.82$ for $\delta_{4S} = 1/5$ and $C_S \approx 10.47$ for
$\delta_{4S} = 1/4$.
\end{theorem}
It is interesting to note that for $S$ obeying the condition of the
theorem, the reconstruction from noiseless data is exact. It is quite
possible that for some matrices $A$, this condition tolerates larger
values of $S$ than $\eqref{eq:d3s}$.

We would like to offer two comments. First, the matrix $A$ is
rectangular with many more columns than rows.  As such, most of its
singular values are zero.  As emphasized earlier, the fact that the
severely ill-posed matrix inversion keeps the perturbation from
``blowing up'' is rather remarkable and perhaps unexpected.

Second, no recovery method can perform fundamentally better for
arbitrary perturbations of size $\epsilon$. To see why this is true,
suppose one had available an {\em oracle} letting us know, in advance,
the support $T_0$ of $x_0$. With this additional information, the
problem is well-posed and one could reconstruct $x_0$ by the method of
Least-Squares for example,
\[
\hat{x} = \begin{cases} (A_{T_0}^*A_{T_0})^{-1} A_{T_0}^* y & \mbox{on } T_0\\
0 & \mbox{elsewhere.} \end{cases}
\]
In the absence of any other information, one could easily argue that
no method would exhibit a fundamentally better performance. Now of
course, $\hat x - x_0 = 0$ on the complement of $T_0$ while on $T_0$
\[
\hat x - x_0 = (A_{T_0}^*A_{T_0})^{-1} A_{T_0}^* e,
\]
and since by hypothesis, the eigenvalues of $A_{T_0}^* A_{T_0}$ are
well-behaved\footnote{Observe the role played by the singular values
  of $A_{T_0}$ in the analysis of the oracle error.}
\[ 
\|\hat{x}-x_0\|_{\ell_2} \approx \|A_{T_0}^* e\|_{\ell_2} \approx
\epsilon, 
\] 
at least for perturbations concentrated in the row space of
$A_{T_0}$. In short, obtaining a reconstruction with an error term
whose size is guaranteed to be proportional to the noise level is the
best one can hope for.

Remarkably, not only can we recover sparse input vectors but one can
also stably recover approximately sparse vectors, as we have the
following companion theorem.  
\begin{theorem}
\label{th:stable2}
Suppose that $x_0$ is an arbitrary vector in $\R^\ncol$ and let $x_{0,S}$
be the truncated vector corresponding to the $S$ largest values of
$x_0$ (in absolute value). Under the hypothesis of Theorem
\ref{th:stable}, the solution $x^\sharp$ to $(P_2)$ obeys
\begin{equation}
\label{eq:lqerror2}
\|x^\sharp - x_0\|_{\ell_2} \leq C_{1,S} 
\cdot \epsilon + C_{2,S} \cdot 
\frac{\|x_0 - x_{0,S}\|_{\ell_1}}{\sqrt{S}}.
\end{equation}
For reasonable values of $\delta_{4S}$ the constants in
\eqref{eq:lqerror} are well behaved; e.g. $C_{1,S} \approx 12.04$ and
$C_{1,S} \approx 8.77$ for $\delta_{4S} = 1/5$.
\end{theorem}
Roughly speaking, the theorem says that minimizing $\ell_1$ stably
recovers the $S$-largest entries of an $\ncol$-dimensional unknown
vector $x$ from $\nrow$ measurements only. 

We now specialize this result to a commonly discussed model in
mathematical signal processing, namely, the class of {\em
  compressible} signals. We say that $x_0$ is compressible if its
entries obey a power law
\begin{equation}
\label{eq:powerlaw}
|x_0|_{(k)} \leq C_r\cdot k^{-r}, 
\end{equation}
where $|x_0|_{(k)}$ is the $k$th largest value of $x_0$ ($|x_0|_{(1)}
\geq |x_0|_{(2)} \geq \ldots \geq |x_0|_{(\ncol)}$), $r > 1$, and
$C_r$ is a constant which depends only on $r$.  Such a model is
appropriate for the wavelet coefficients of a piecewise smooth signal,
for example.  If $x_0$ obeys \eqref{eq:powerlaw}, then
\[
\frac{\|x_0 - x_{0,S}\|_{\ell_1}}{\sqrt{S}} \le C'_r \cdot S^{-r+1/2}.
\]
Observe now that in this case
\[
\|x_0 - x_{0,S}\|_{\ell_2} \le C''_r \cdot S^{-r+1/2},
\]
and for generic elements obeying \eqref{eq:powerlaw}, there are no
fundamentally better estimates available. Hence, we see that with
$\nrow$ measurements only, we achieve an approximation error which is
almost as good as that one would obtain by knowing everything about
the signal $x_0$ and selecting its $S$-largest entries. 

As a last remark, we would like to point out that in the noiseless
case, Theorem \ref{th:stable2} improves upon an earlier result from
Cand\`es and Tao, see also \cite{donoho04co}; it is sharper in the
sense that 1) this is a deterministic statement and there is no
probability of failure, 2) it is universal in that it holds for all
signals, 3) it gives upper estimates with better bounds and constants,
and 4) it holds for a wider range of values of $S$.

%----------------------------------------------------------------------
\subsection{Examples}
\label{sec:examples}

It is of course of interest to know which matrices obey the uniform
uncertainty principle with good isometry constants.  Using tools from
random matrix theory, \cite{candes04ro,candes04ne,donoho04fo} give
several examples of matrices such that \eqref{eq:d3s} holds for $S$ on
the order of $\nrow$ to within log factors.  Examples include (proofs
and additional discussion can be found in \cite{candes04ne}):

\begin{itemize}
\item {\em Random matrices with i.i.d.~entries}. Suppose the entries
  of $A$ are i.i.d.~Gaussian with mean zero and variance $1/\nrow$,
  then \cite{szarek91co,candes04ne,donoho04fo} show that the condition
  for Theorem~\ref{th:stable} holds with overwhelming probability when
\[
S\leq C\cdot \nrow/\log(\ncol/\nrow).
\]
In fact, \cite{candes04de} gives numerical values for the constant $C$
as a function of the ratio $\nrow/\ncol$. The same conclusion applies
to binary matrices with independent entries taking values $\pm
1/\sqrt{\nrow}$ with equal probability.

% Since the Gaussian distribution is completely isotropic, choosing a
% $\nrow\times \ncol$ matrix in this manner is in some sense equivalent to
% choosing a $\nrow$ dimensional subspace onto which project the sparse
% signal $x_0$.  That our matrix is properly conditioned with
% overwhelming probability means that we can successfully recover $x_0$
% from its projection onto almost any $\nrow$ dimensional subspace.

% \item {\em Binary ensemble.]  Now let the entries of $A$ be binary, taking
%   values $\pm 1/\sqrt{\nrow}$ with equal probability.  Just as in the
%   Gaussian case, the condition for Theorem~\ref{th:stable} holds with
%   overwhelming probability when $S\leq C\cdot \nrow/\log (\ncol/\nrow)$.
  
\item {\em Fourier ensemble.} Suppose now that $A$ is obtained by
  selecting $\nrow$ rows from the $\ncol \times \ncol$ discrete
  Fourier transform and renormalizing the columns so that they are
  unit-normed.  If the rows are selected at random, the condition for
  Theorem~\ref{th:stable} holds with overwhelming probability for
  $S\leq C\cdot \nrow/(\log \ncol)^6$ \cite{candes04ne}. (For
  simplicity, we have assumed that $A$ takes on real-valued entries
  although our theory clearly accommodates complex-valued matrices so
  that our discussion holds for both complex and real-valued Fourier
  transforms.)

% The Fourier
%   ensemble has complex entries.  A $\nrow\times\ncol$ {\em real
%     Fourier ensemble} can be created by separating the real and
%   imaginary (cosine and sine) parts of a standard $\nrow/2\times\ncol$
%   Fourier ensemble.
    
  This case is of special interest as reconstructing a digital signal
  or image from incomplete Fourier data is an important inverse
  problem with applications in biomedical imaging (MRI and
  tomography), Astrophysics (interferometric imaging), and geophysical
  exploration.
  
\item {\em General orthogonal measurement ensembles.}  Suppose $A$ is
  obtained by selecting $\nrow$ rows from an $\ncol$ by $\ncol$
  orthonormal matrix $U$ and renormalizing the columns so that they
  are unit-normed.  Then \cite{candes04ne} shows that if the rows are
  selected at random, the condition for Theorem~\ref{th:stable} holds
  with overwhelming probability provided 
  \begin{equation}
    \label{eq:cond-incoherence}
    S\leq C \cdot \frac{1}{\mu^2} \cdot \frac{\nrow}{(\log \ncol)^6},
  \end{equation}
  where $\mu := \sqrt{m} \, \max_{i,j} |U_{i,j}|$. Observe that for
  the Fourier matrix, $\mu = 1$, and thus
  \eqref{eq:cond-incoherence} is an extension of the Fourier ensemble.
 
  This fact is of significant practical relevance because in many
  situations, signals of interest may not be sparse in the time domain
  but rather may be (approximately) decomposed as a sparse
  superposition of waveforms in a fixed orthonormal basis $\Psi$;
  e.g. in a nice wavelet basis. Suppose that we use as test signals a
  set of $\nrow$ vectors taken from a second orthonormal basis
  $\Phi$. We then solve $(P_1)$ in the coefficient domain
\[
(P_1^\prime)\quad\quad \min~\|\alpha\|_{\ell_1} \quad\text{subject to}
\quad A \alpha = y, 
\]
where $A$ is obtained by extracting $\nrow$ rows from the orthonormal
matrix $U = \Phi \Psi^*$. The recovery condition then depends on the
{\em mutual coherence} $\mu$ between the measurement basis $\Phi$ and
the sparsity basis $\Psi$ which measures the similarity between $\Phi$
and $\Psi$; $\mu(\Phi,\Psi) = \sqrt{m}\, \max~|\<\phi_k,\psi_j\>|$,
$\phi_k\in\Phi$, $\psi_j\in\Psi$.
% the restricted isometry constants of the matrix $\Phi_\Omega\Psi^*$.
% The $\delta_S$, in turn, will depend on the {\em incoherence} $\mu$
% between the measurement basis $\Phi$ and the sparsity basis $\Psi$.
% \[
% \mu = \max~|\<\phi_k,\psi_j\>|\quad \phi_k\in\Phi, \psi_j\in\Psi.
% \]
% The parameter $\mu$ is a basic measure of the similarity of the basis
% $\Phi$ to the basis $\Psi$.  In \cite{candes05sp}, it is shown that if
% we select the subset of measurement vectors $\Omega$ at random, then
% with overwhelming probability the condition for
% Theorem~\ref{th:stable} holds for $S\leq C\cdot (\ncol/\mu^2)\cdot \nrow/(\log
% \ncol)^3$.
\end{itemize}

%----------------------------------------------------------------------
\subsection{Prior work and innovations}

The problem of recovering a sparse vector by minimizing $\ell_1$ under
linear equality constraints has recently received much attention,
mostly in the context of {\em Basis Pursuit}, where the goal is to
uncover sparse signal decompositions in overcomplete dictionaries. We
refer the reader to \cite{donoho01un,donoho03op} and the references
therein for a full discussion. 

We would especially like to note two works by Donoho, Elad, and
Temlyakov \cite{donoho04st}, and Tropp \cite{tropp05ju} that also
study the recovery of sparse signals from noisy observations by
solving $(P_2)$ (and other closely related optimization programs), and
give conditions for stable recovery.  In \cite{donoho04st}, the
sparsity constraint on the underlying signal $x_0$ depends on the
magnitude of the maximum entry of the Gram matrix $M(A) = \max_{i,j :
  i\not=j} |(A^*A)|_{i,j}$.  Stable recovery occurs when the number of
nonzeros is at most $(M^{-1} + 1)/4$.  For instance, when $A$ is a
Fourier ensemble and $\nrow$ is on the order of $\ncol$, we will have
$M$ at least of the order $1/\sqrt{\nrow}$ (with high probability),
meaning that stable recovery is known to occur when the number of
nonzeros is about at most $O(\sqrt{\nrow})$.  In contrast, the
condition for Theorem~\ref{th:stable} will hold when this number is
about $\nrow/(\log \ncol)^6$, due to the range of support sizes for
which the uniform uncertainty principle holds.  In \cite{tropp05ju}, a
more general condition for stable recovery is derived.  For the
measurement ensembles listed in the previous section, however, the
sparsity required is still on the order of $\sqrt{\nrow}$ in the
situation where $\nrow$ is comparable to $\ncol$. In other words, {\em
  whereas these results require at least $O(\sqrt{\ncol})$
  observations per unknown, our results show that---ignoring log-like
  factors---only $O(1)$ are, in general, sufficient.}

More closely related is the very recent work of Donoho
\cite{donoho04fo2} who shows a version of \eqref{eq:lqerror} in the
case where $A \in \R^{\nrow \times \ncol}$ is a Gaussian matrix with
$\nrow$ proportional to $\ncol$, with unspecified constants for both
the support size and that appearing in \eqref{eq:lqerror}.  Our main
claim is on a very different level since it is (1) deterministic (it
can of course be specialized to random matrices), and (2) widely
applicable since it extends to any matrix obeying the condition
$\delta_{3S} + 3\delta_{4S} < 2$. In addition, the argument underlying
Theorem \ref{th:stable} is short and simple, giving precise and
sharper numerical values. Finally, we would like to point out
connections with fascinating ongoing work which develops fast
randomized algorithms for sparse Fourier transforms
\cite{gilbert04be,zou04th}. Suppose $x_0$ is a fixed vector with
$|T_0|$ nonzero terms, for example. Then \cite{gilbert04be} shows that
it is possible to randomly sample the frequency domain $|T_0|
\text{poly}(\log m)$ times ($\text{poly}(\log m)$ denotes a polynomial
term in $\log m$), and reconstruct $x_0$ from these frequency data
with positive probability. We do not know whether these algorithms are
stable in the sense described in this paper, and whether they can be
modified to be universal, i.e. reconstruct all signals of small
support.

%----------------------------------------------------------------------
\section{Proofs}
\label{sec:proof}

\subsection{Proof of Theorem~\ref{th:stable}: sparse case}

\begin{figure}
  \centering 
\includegraphics[scale = .55]{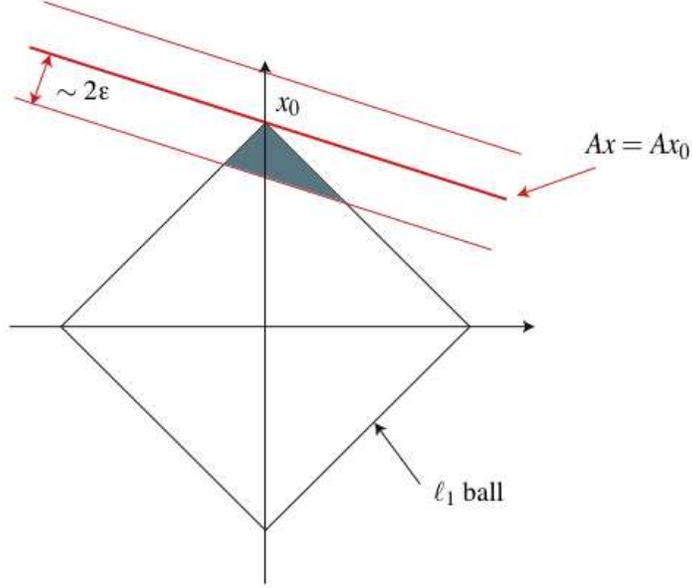}
\caption{\small\sl Geometry in $\R^2$.  Here, the point $x_0$ is a
  vertex of the $\ell_1$ ball and the shaded area represents the set
  of points obeying both the tube and the cone constraints. By showing
  that every vector in the cone of descent at $x_0$ is approximately
  orthogonal to the nullspace of $A$, we will ensure that $x^\sharp$
  is not too far from $x_0$.}
\label{fig:geom}
\end{figure}

The proof of the theorem makes use of two geometrical special facts
about the solution $x^\sharp$ to $(P_2)$.
\begin{enumerate}
\item {\em Tube constraint.}  First, $Ax^\sharp$ is within
  $2\epsilon$ of the ``noise free'' observations $Ax_0$ thanks to the
  triangle inequality
\begin{equation}
\label{eq:tube}
\|Ax^\sharp - Ax_0\|_{\ell_2} \leq 
\|Ax^\sharp - y\|_{\ell_2} + \|Ax_0 - y\|_{\ell_2} \leq 2\epsilon.
\end{equation}
Geometrically, this says that $x^\sharp$ is known to be in a cylinder
around the $\nrow$-dimensional plane $A x_0$.
\item {\em Cone constraint.}  Since $x_0$ is feasible, we must have
  $\|x^\sharp\|_{\ell_1} \leq \|x_0\|_{\ell_1}$.  Decompose $x^\sharp$
  as $x^\sharp = x_0 + h$.  As observed in \cite{donoho01un}
\[
\|x_0\|_{\ell_1} - \|h_{T_0}\|_{\ell_1} + \|h_{T^c_0}\|_{\ell_1} \leq
\|x_0 + h\|_{\ell_1} \leq \|x_0\|_{\ell_1}, 
\]
where $T_0$ is the support of $x_0$, and $h_{T_0}(t) = h(t)$ for $t\in
T_0$ and zero elsewhere (similarly for $h_{T^c_0}$). Hence, $h$ obeys
the cone constraint
\begin{equation}
\label{eq:conecond}
\|h_{T^c_0}\|_{\ell_1} \leq  \|h_{T_0}\|_{\ell_1}
\end{equation}
which expresses the geometric idea that $h$ must lie in the cone of
descent of the $\ell_1$-norm at $x_0$.
\end{enumerate}
Figure~\ref{fig:geom} illustrates both these geometrical constraints.
Stability follows from the fact that the intersection between
\eqref{eq:tube} ($\|A h \|_{\ell_2} \le 2 \epsilon$) and
\eqref{eq:conecond} is a set with small radius. The reason why this
holds is because every vector $h$ in the $\ell_1$-cone
\eqref{eq:conecond} is approximately orthogonal to the nullspace of
$A$. We shall prove that $\|Ah\|_{\ell_2} \approx \|h\|_{\ell_2}$ and
together with \eqref{eq:tube}, this establishes the theorem.

%-----------------
%-----------------

We begin by dividing $T^c_0$ into subsets of size $M$ (we will choose
$M$ later) and enumerate $T^c_0$ as $n_1,n_2,\ldots,n_{N-|T_0|}$ in
decreasing order of magnitude of $h_{T^c_0}$.  Set $T_j = \{n_\ell,
(j-1)M+1 \leq \ell \leq jM\}$.  That is, $T_1$ contains the indices
of the $M$ largest coefficients of $h_{T^c_0}$, $T_2$ contains the
indices of the next $M$ largest coefficients, and so on.

With this decomposition, the $\ell_2$-norm of $h$ is concentrated on
$T_{01} = T_0\cup T_1$.  Indeed, the $k$th largest value of
$h_{T^c_0}$ obeys
\[
|h_{T^c_0}|_{(k)} ~\leq ~\|h_{T^c_0}\|_{\ell_1}/k
\]
and, therefore, 
\[
\|h_{T^c_{01}}\|^2_{\ell_2} ~\leq ~ \|h_{T^c_0}\|^2_{\ell_1}
\sum_{k=M+1}^\ncol 1/k^2 ~\leq ~\|h_{T^c_0}\|^2_{\ell_1}/M.
\]
Further, the $\ell_1$-cone constraint gives
\[
\|h_{T^c_{01}}\|^2_{\ell_2} ~\leq ~\|h_{T_0}\|^2_{\ell_1}/M ~\leq
~\|h_{T_0}\|^2_{\ell_2}\cdot|T_0|/M
\]
and thus
\begin{equation}
\label{eq:hnorm}
\|h\|^2_{\ell_2} ~= ~\|h_{T_{01}}\|^2_{\ell_2} + 
\|h_{T^c_{01}}\|^2_{\ell_2} ~\leq ~ 
(1+|T_0|/M)\cdot \|h_{T_{01}}\|^2_{\ell_2}.
\end{equation}

Observe now that 
\begin{eqnarray*}
\|Ah\|_{\ell_2}  = \|A_{T_{01}}h_{T_{01}} + \sum_{j\geq 2} A_{T_j}h_{T_j} \|_{\ell_2} 
 & \ge & 
\|A_{T_{01}}h_{T_{01}}\|_{\ell_2} - \|\sum_{j\geq 2} A_{T_j}h_{T_j} \|_{\ell_2} \\
 & \geq & \|A_{T_{01}}h_{T_{01}}\|_{\ell_2} - \sum_{j\geq 2} \|A_{T_j}h_{T_j}\|_{\ell_2} \\
 & \geq & \sqrt{1-\delta_{M+|T_0|}}\,\|h_{T_{01}}\|_{\ell_2} - 
 \sqrt{1+\delta_M}\, \sum_{j\geq 2} \|h_{T_j}\|_{\ell_2}. 
\end{eqnarray*}
Set $\rho = |T_0|/M$. As we shall see later,
\begin{equation}
\label{eq:sumnorm}
\sum_{j\geq 2} \|h_{T_j}\|_{\ell_2} \le  \sqrt{\rho}\cdot \|h_{T_0}\|_{\ell_2}  
\end{equation}
which gives
\begin{equation}
  \label{eq:Ct0}
\|Ah\|_{\ell_2} \ge C_{|T_0|,M} \cdot \|h_{T_{01}}\|_{\ell_2}, \quad 
C_{|T_0|,M} := \sqrt{1-\delta_{M+|T_0|}} -
\sqrt{\rho}\, \sqrt{1+\delta_M}.
\end{equation}
It then follows from \eqref{eq:hnorm} and $\|Ah\|_{\ell_2} \leq
2\epsilon$ that
\begin{equation}
\label{eq:good}
\|h\|_{\ell_2} \le \sqrt{1+ \rho} \cdot \|h_{T_{01}}\|_{\ell_2} \le
\frac{\sqrt{1+ \rho}}{C_{|T_0|, M}} \cdot \|A h\|_{\ell_2} \le 
\frac{2\sqrt{1+ \rho}}{C_{|T_0|, M}} \cdot \epsilon,   
\end{equation}
provided that the denominator is of course positive.
 
We may specialize \eqref{eq:good} and take $M = 3|T_0|$. The
denominator is positive if $\delta_{3|T_0|} + 3 \delta_{4|T_0|} < 2$ (this is
true if $\delta_{4|T_0|} < 1/2$, say) which proves the theorem.  Note
that if $\delta_{4S}$ is a little smaller, the constant in
\eqref{eq:lqerror} is not large. For $\delta_{4S}\leq 1/5$, $C_S
\approx 8.82$, while for $\delta_{4S}\leq 1/4$, $C_S \approx 10.47$ as
claimed.

It remains to argue about \eqref{eq:sumnorm}.  Observe that by
construction, the magnitude of each coefficient in $T_{j+1}$ is less
than the average of the magnitudes in $T_j$:
\[
|h_{T_{j+1}}(t)| ~\leq \|h_{T_j}\|_{\ell_1}/M.
\]
Then
\[
\|h_{T_{j+1}}\|^2_{\ell_2} ~ \leq ~ \|h_{T_j}\|^2_{\ell_1}/M
\]
and \eqref{eq:sumnorm} follows from
\[
\sum_{j\geq 2} \|h_{T_j}\|_{\ell_2} ~ \leq ~ \sum_{j\geq 1}
\|h_{T_j}\|_{\ell_1}/\sqrt{M} ~\leq ~ \|h_{T_0}\|_{\ell_1}/\sqrt{M}
\leq \sqrt{|T_0|/M}\cdot\|h_{T_0}\|_{\ell_2}.
\]

\subsection{Proof of Theorem~\ref{th:stable2}: general case}

Suppose now that $x_0$ is arbitrary. We let $T_0$ be the indices of
the largest $|T_0|$ coefficients of $x_0$ (the value $|T_0|$ will be
decided later) and just as before, we divide up $T_0^c$ into sets
$T_1, \ldots, T_J$ of equal size $|T_j| = M$, $j \ge 1$, by
decreasing order of magnitude. The cone constraint \eqref{eq:conecond}
may not hold but a variation does. Indeed, $x = x_0 + h$ is
feasible and, therefore,
\[ 
\|x_{0,T_0}\|_{\ell_1} - \|h_{T_0}\|_{\ell_1} -
\|x_{0,T_0^c}\|_{\ell_1} + \|h_{T_0^c}\|_{\ell_1} \le \|x_{0,T_0} +
h_{T_0}\|_{\ell_1} + \|x_{0,T_0^c} + h_{T_0^c}\|_{\ell_1} \le
\|x_0\|_{\ell_1},
\]
which gives
\begin{equation}
\label{eq:conecond2}
\|h_{T^c_0}\|_{\ell_1} \leq \|h_{T_0}\|_{\ell_1} + 2
\|x_{0,T_0^c}\|_{\ell_1}. 
\end{equation}

The rest of the argument now proceeds essentially as before. First,
$h$ is in the some sense concentrated on $T_{01} = T_0 \cup T_1$ since
with the same notations
\[ 
\|h_{T^c_{01}}\|_{\ell_2} \le \frac{\|h_{T_0}\|_{\ell_1} + 2
\|x_{0,T_0^c}\|_{\ell_1}}{\sqrt{M}} \le \sqrt{\rho} \cdot 
\left(\|h_{T_0}\|_{\ell_2} + 
\frac{2 \|x_{0,T_0^c}\|_{\ell_1}}{\sqrt{|T_0|}}\right), 
\]
which in turn implies 
\begin{equation}
  \label{eq:h2}
\|h\|_{\ell_2} \le (1 + \sqrt{\rho}) \|h_{T_{01}}\|_{\ell_2}
+ 2 \sqrt{\rho} \cdot \eta, \quad \eta :=
\|x_{0,T_0^c}\|_{\ell_1}/\sqrt{|T_0|}.
\end{equation}
Better estimates via Pythagoras' formula are of course possible (see
\eqref{eq:hnorm}) but we ignore such refinements in order to keep the
argument as simple as possible.  Second, the same reasoning as before
gives
$$ \sum_{j \ge 2} \|h_{T_j}\|_{\ell_2} \le
\frac{\|h_{T_{01}^c}\|_{\ell_1}}{\sqrt{M}} \le \sqrt{\rho} \cdot
(\|h_{T_0}\|_{\ell_2} + 2 \eta)
$$
and thus
\[ 
\|Ah\|_{\ell_2} \ge C_{|T_0|,M} \cdot \|h_{T_{01}}\|_{\ell_2} -
2\sqrt{\rho}\, \sqrt{1+\delta_M} \cdot \eta, 
\] 
where $C_{|T_0|,M}$ is the same as in \eqref{eq:Ct0}.  Since $\|A h\|
\le 2\epsilon$, we again conclude that
\[ 
\|h_{T_{01}}\|_{\ell_2} \le \frac{2}{C_{|T_0|,M}} \cdot (\epsilon +
\sqrt{\rho}\, \sqrt{1+\delta_M} \, \eta),
\]
(note that the constant in front of the $\epsilon$ factor is the same
as in the truly sparse case) and the claim \eqref{eq:lqerror2} follows
from \eqref{eq:h2}. Specializing the bound to $M = 3 |T_0|$ and
assuming that $\delta_S \le 1/5$ gives the numerical values reported
in the statement of the theorem.

%----------------------------------------------------------------------
\section{Numerical Examples}
\label{sec:numerical}

This section illustrates the effectiveness of the recovery by means of
a few simple numerical experiments.  Our simulations demonstrate that
in practice, the constants in \eqref{eq:lqerror} and
\eqref{eq:lqerror2} seem to be quite low.

Our first series of experiments is summarized in
Tables~\ref{tab:sparse} and \ref{tab:compress}.  In each experiment, a
length $1024$ signal was measured with the same $300\times 1024$
Gaussian measurement ensemble.  The measurements were then corrupted
by additive white Gaussian noise: $y_k = \<x_0,a_k\> + e_k$ with
$e_k\sim \mathcal{N}(0,\sigma^2)$ for various noise levels $\sigma$.
The squared norm of the error $\|e\|^2_{\ell_2}$ is a chi-square
random variable with mean $\sigma^2 \nrow$ and standard deviation
$\sigma^2 \sqrt{2\nrow}$; owing to well known concentration
inequalities, the probability that $\|e\|^2_{\ell_2}$ exceeds its mean
plus two or three standard deviations is small. We then solve $(P_2)$
with
\begin{equation}
\label{eq:epsilon}
\epsilon^2 = \sigma^2(\nrow + \lambda \sqrt{2\nrow}) 
\end{equation}
and select $\lambda = 2$ although other choices are of course
possible.

\begin{figure}
\centerline{
\begin{tabular}{cc}
\includegraphics[width=2.5in]{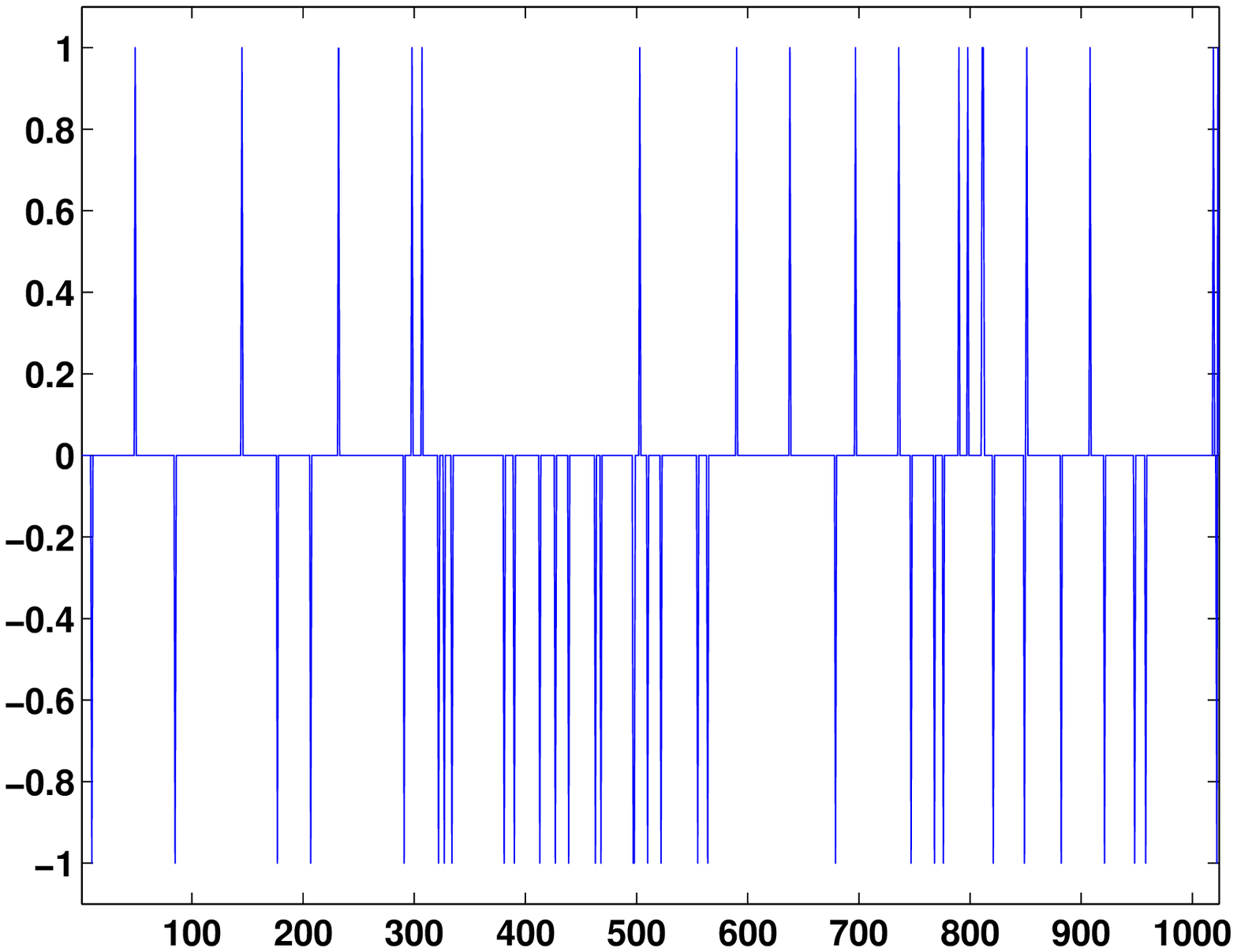} & 
\includegraphics[width=2.5in]{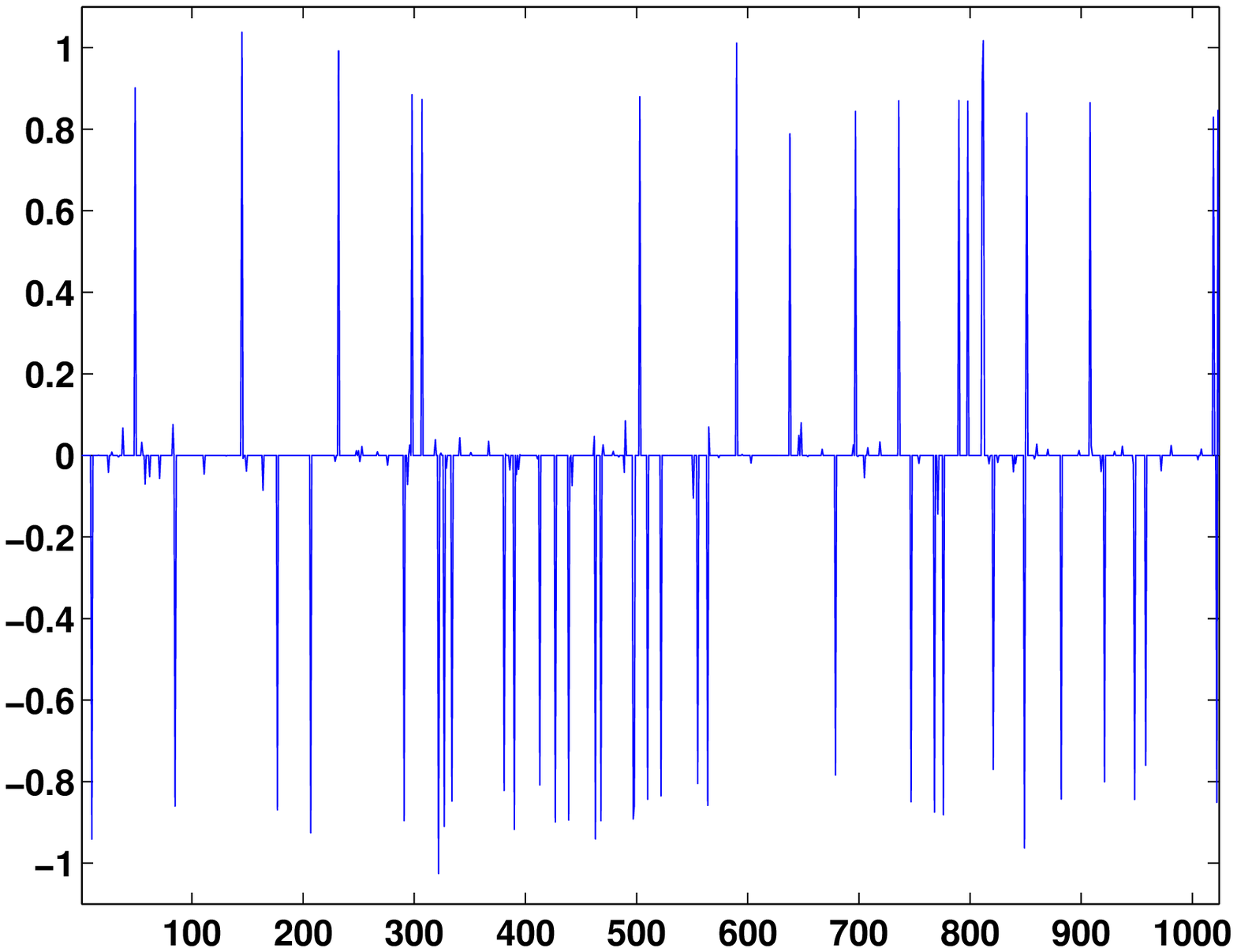} \\
(a) & (b) \\
\includegraphics[width=2.5in]{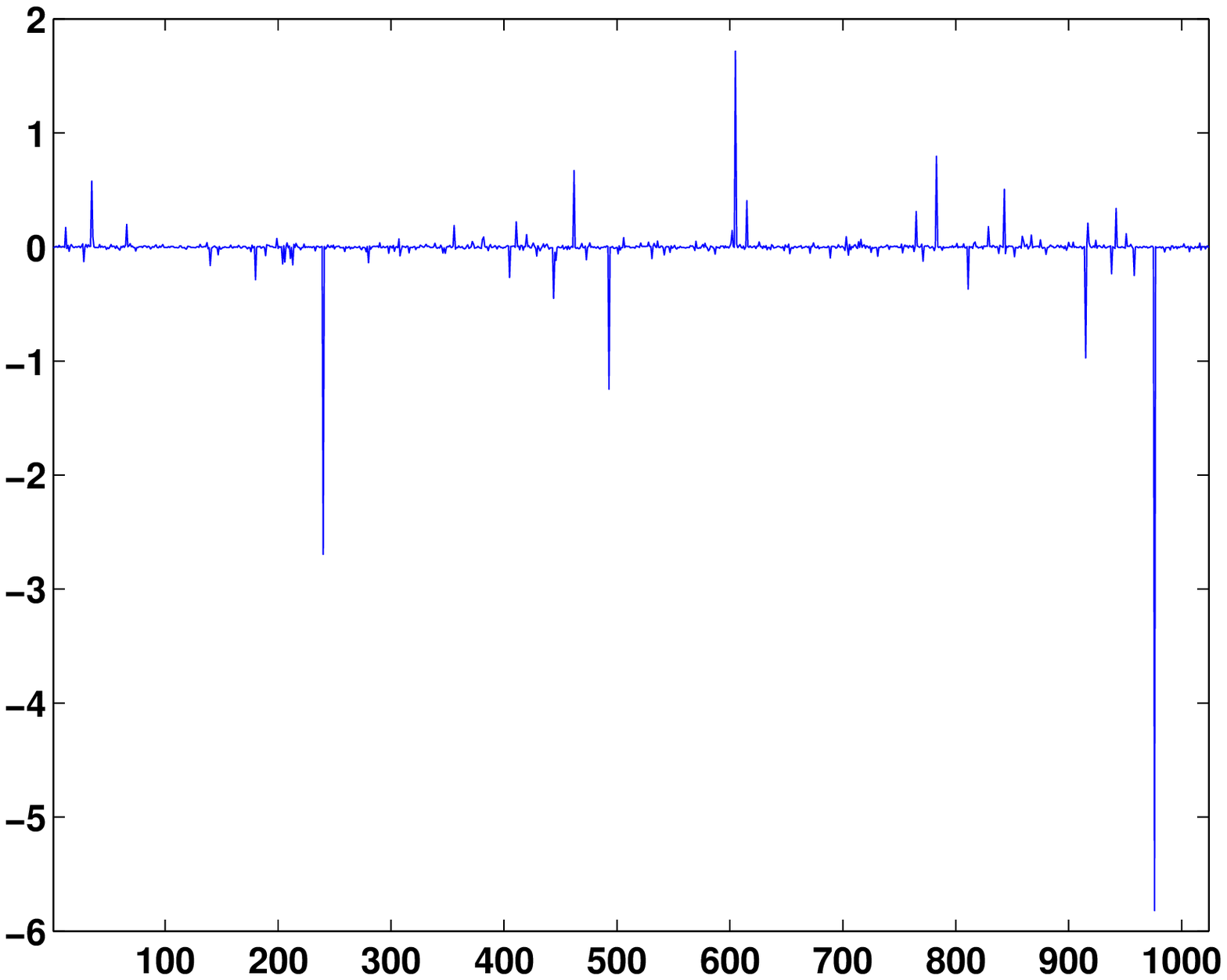} & 
\includegraphics[width=2.5in]{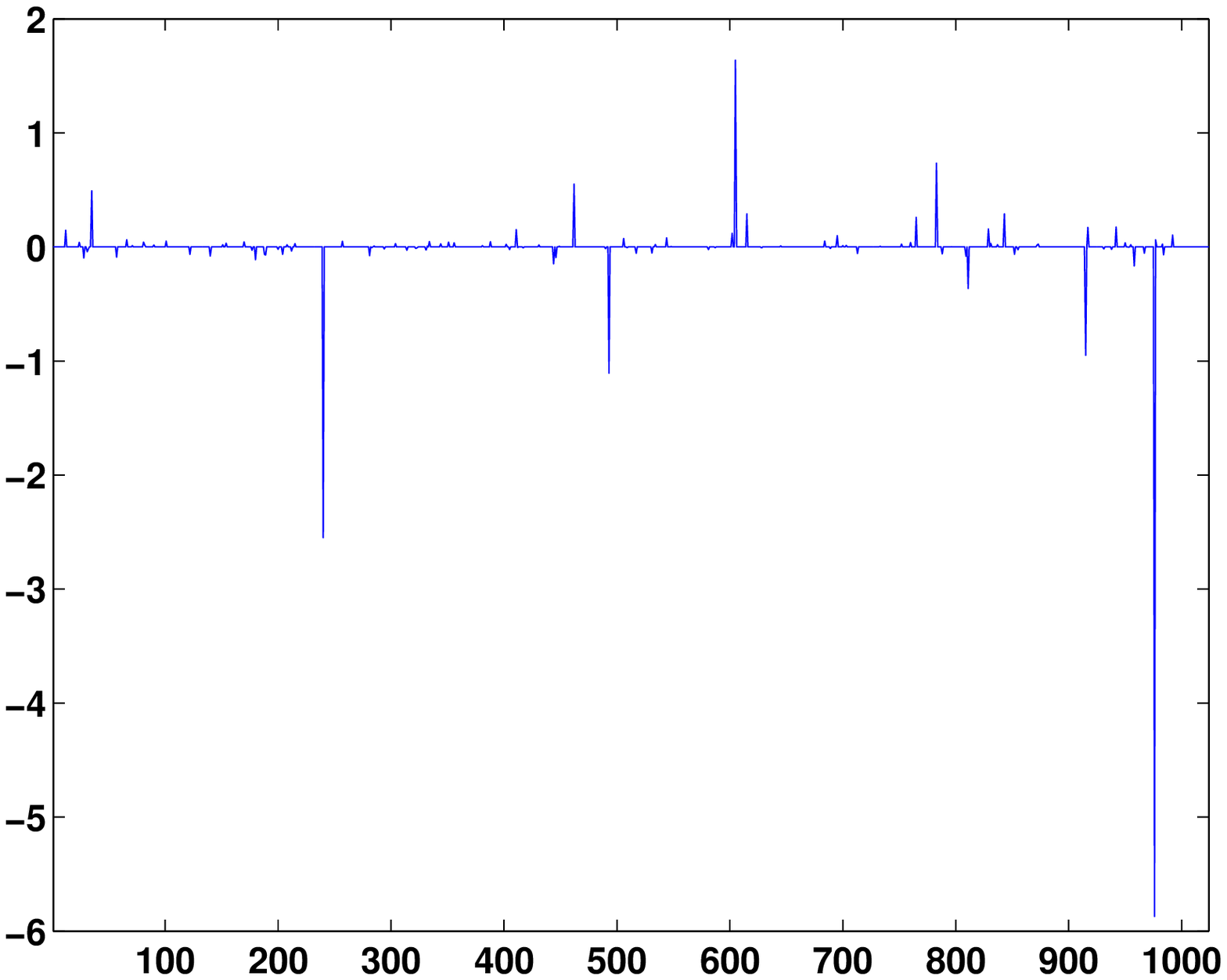} \\
(c) & (d)
\end{tabular}
}
\caption{\small\sl (a) Example of a sparse signal used in the 1D experiments.  There are $50$ non-zero coefficients taking values $\pm 1$.
(b) Sparse signal recovered from noisy measurements with $\sigma = 0.05$.
(c) Example of a compressible signal used in the 1D experiments.
(d) Compressible signal recovered from noisy measurements with $\sigma=0.05$.
}
\label{fig:pics}
\end{figure}

Table~\ref{tab:sparse} charts the results for sparse signals with $50$
nonzero components.  Ten signals were generated by choosing $50$
indices uniformly at random, and then selecting either $-1$ or $1$ at
each location with equal probability.  An example of such a signal is
shown in Figure~\ref{fig:pics}(a).  Previous experiments
\cite{candes04ro} have demonstrated that we were empirically able to
recover such signals perfectly from $300$ noiseless Gaussian
measurements, which is indeed the case for each of the $10$ signals
considered here.  The average value of the recovery error (taken over
the $10$ signals) is recorded in the bottom row of
Table~\ref{tab:sparse}.  In this situation, the constant in
\eqref{eq:lqerror} appears to be less than $2$.

Table~\ref{tab:compress} charts the results for $10$ compressible signals whose components are all non-zero, but decay as in \eqref{eq:powerlaw}.  The signals were generated by taking a fixed sequence
\begin{equation}
\label{eq:csignals}
x_{\mathrm{sort}}(t) = (5.819)\cdot t^{-10/9},
\end{equation}
randomly permuting it, and multiplying by a random sign sequence (the
constant in \eqref{eq:csignals} was chosen so that the norm of the
compressible signals is the same --- $\sqrt{50}$ --- as the sparse
signals in the previous set of experiments).  An example of such a
signal is shown in Figure~\ref{fig:pics}(c).  Again, $10$ such signals
were generated, and the average recovery error recorded in the bottom
row of Table~\ref{tab:compress}.

For small values of $\sigma$, the recovery error is dominated by the
approximation error --- the second term on the right hand side of
\eqref{eq:lqerror2}.  As a reference, the $50$ term nonlinear
approximation errors of these compressible signals is around $0.47$;
at low signal-to-noise ratios our recovery error is about $1.5$ times
this quantity.  As the noise power gets large, the recovery error
becomes less than $\epsilon$, just as in the sparse case.

\begin{figure}
\centerline{
\begin{tabular}{ccc}
\includegraphics[width=2in]{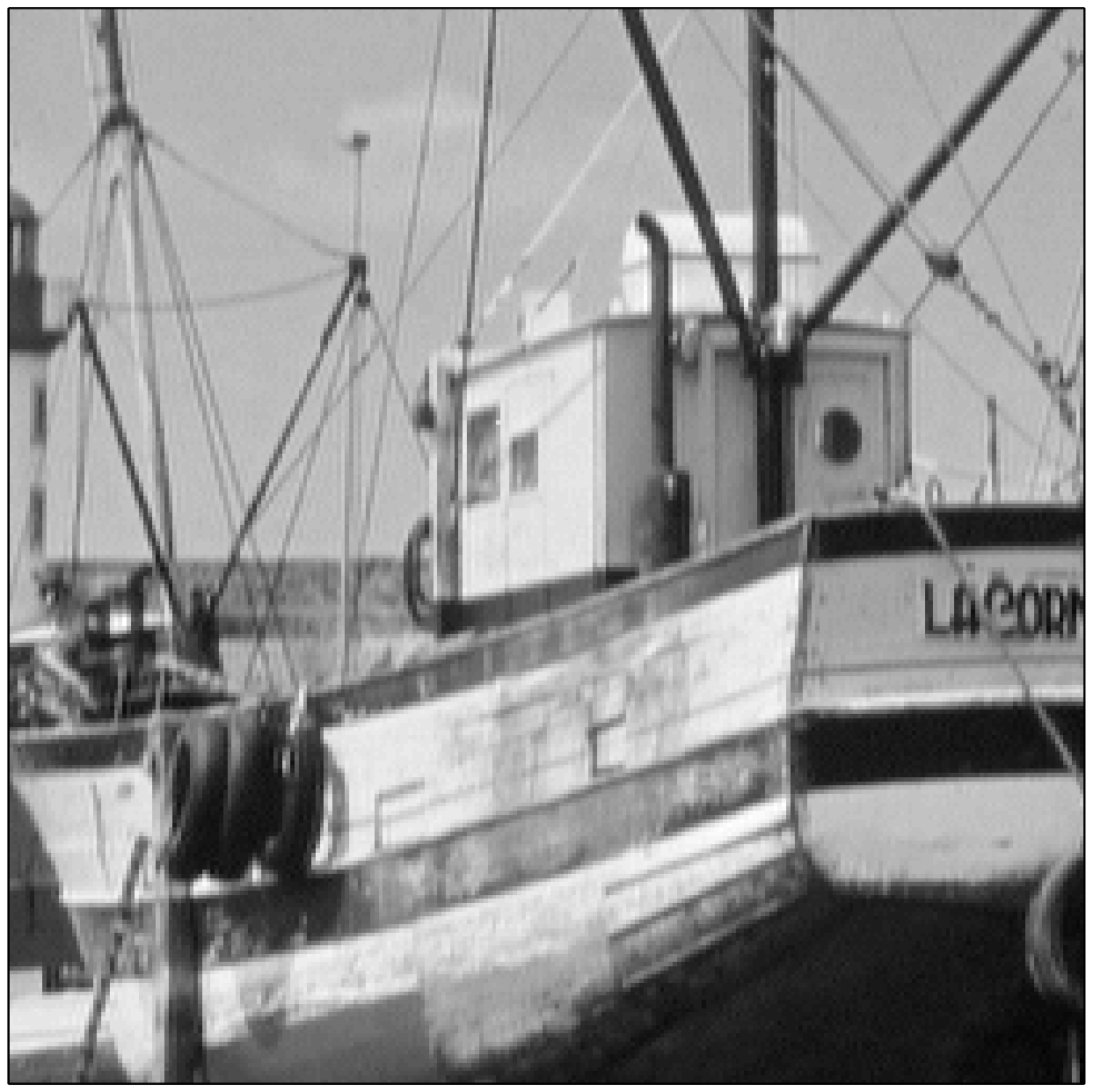} &
\includegraphics[width=2in]{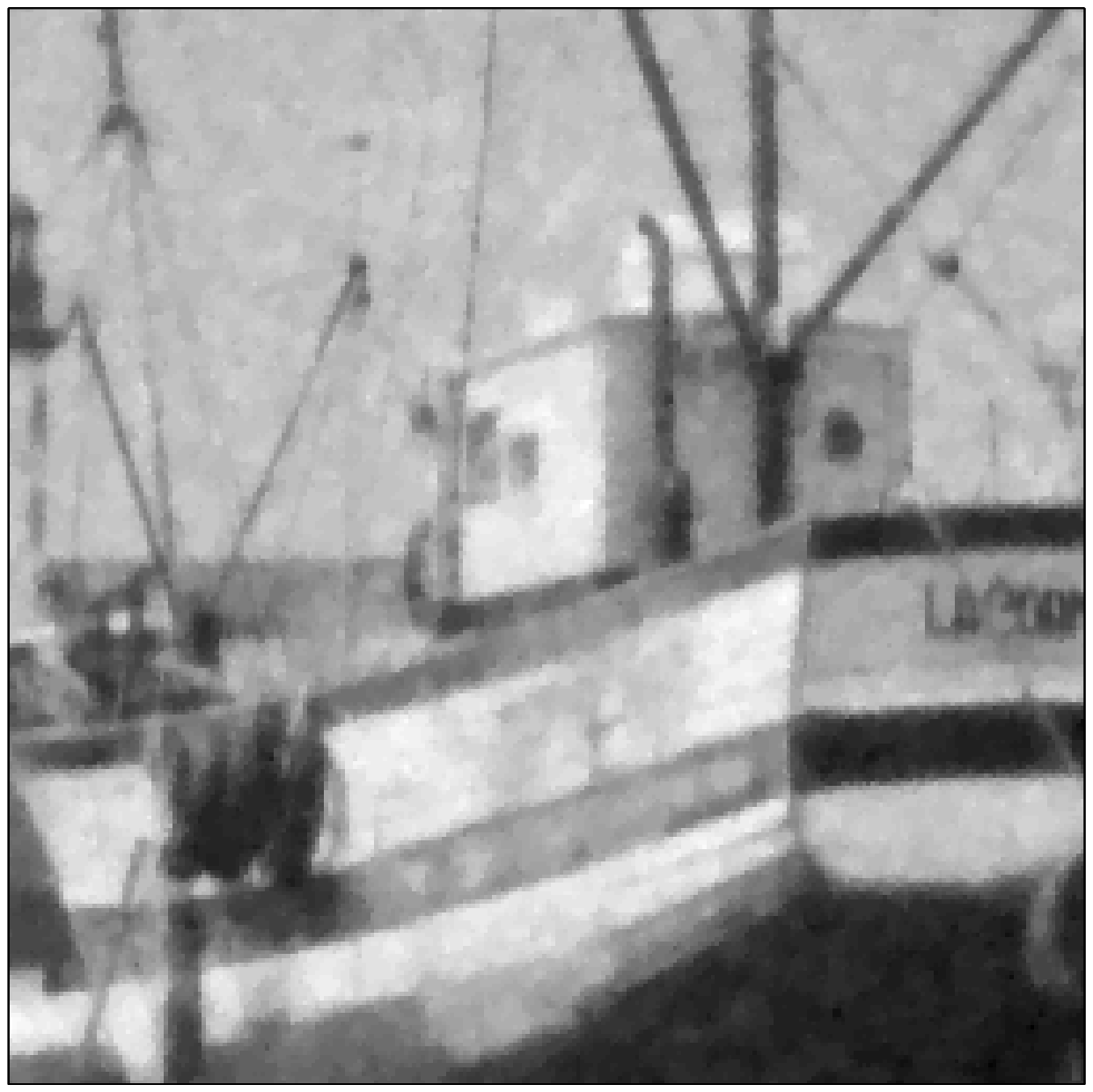} &
\includegraphics[width=2in]{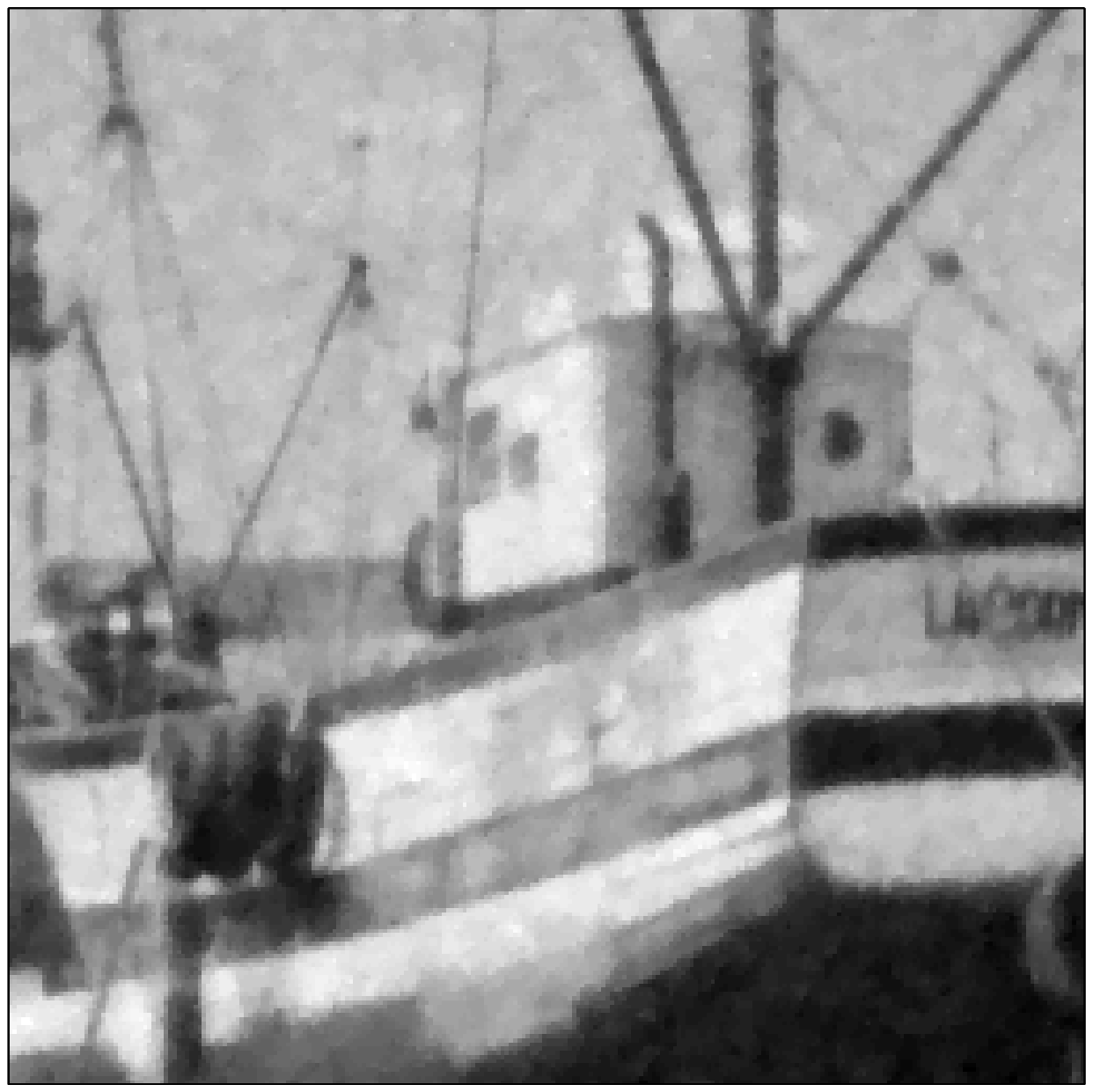} \\
(a) & (b) & (c)
\end{tabular}
}
\caption{\small\sl (a) Original $256\times 256$ {\em Boats} image.  
(b) Recovery via $(TV)$ from $\nrow=25000$ measurements corrupted with Gaussian noise.
(c) Recovery via $(TV)$ from $\nrow=25000$ measurements corrupted by round-off error.
In both cases, the reconstruction error is less than the sum of the nonlinear approximation and measurement errors.
}
\label{fig:image_example}
\end{figure}

Finally, we apply our recovery procedure to realistic imagery.
Photograph-like images, such as the $256\times 256$ pixel {\em Boats}
image shown in Figure~\ref{fig:image_example}(a), have wavelet
coefficient sequences that are compressible (see \cite{devore92im}).
The image is a $65536$ dimensional vector, making the standard
Gaussian ensemble too unwieldy\footnote{Storing a double precision
  $25000\times65536$ matrix would use around $13.1$ gigabytes of
  memory, about the capacity of three standard DVDs.}.  Instead, we
make $25000$ measurements of the image using a {\em scrambled} real
Fourier ensemble; that is, the test functions $a_k(t)$ are real-valued
sines and cosines (with randomly selected frequencies) which are
temporally scrambled by randomly permuting the $\ncol$ time points. In
other words, this ensemble is obtained from the (real-valued) Fourier
ensemble by a random permutation of the columns.  For our purposes here, the
test functions behave like a Gaussian ensemble in the sense that from
$\nrow$ measurements, one can recover signals with about $\nrow/5$
nonzero components exactly from noiseless data. 
There is a computational advantage as well,
since we can apply $A$ and its adjoint
$A^T$ to an arbitrary vector by means of an $m$ point FFT.  To recover
the wavelet coefficients of the object, we simply solve
\[
(P^\prime_2)\qquad \min~\|\alpha\|_{\ell_1} \quad \text{subject~to}
\quad \|AW^*\alpha - y\|_{\ell_2} \leq \epsilon,
\]
where $A$ is the scrambled Fourier ensemble, and $W$ is the discrete
Daubechies-8 orthogonal wavelet transform.

\begin{figure}
\centerline{
\begin{tabular}{ccc}
\includegraphics[width=2in]{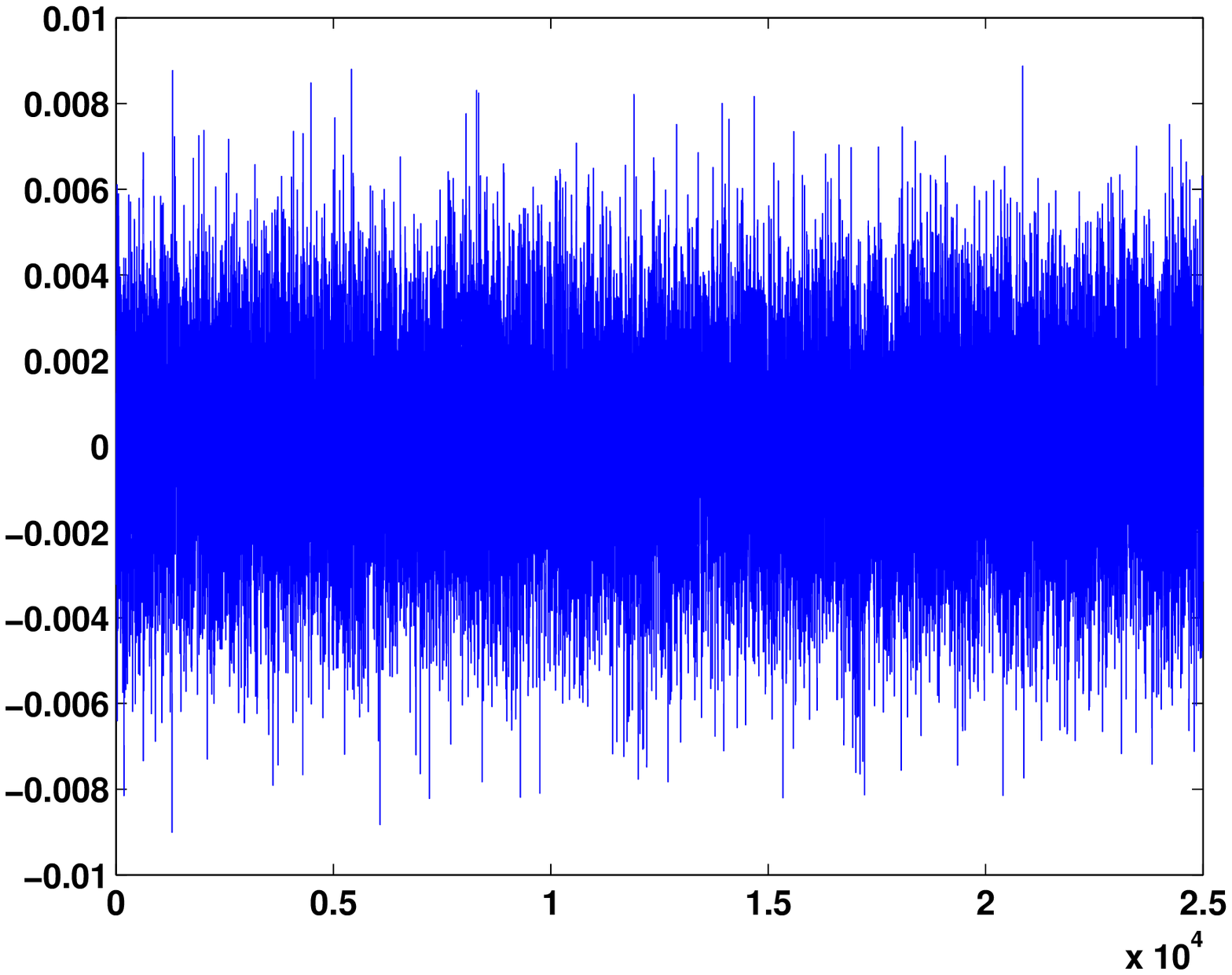} &
\includegraphics[width=2in]{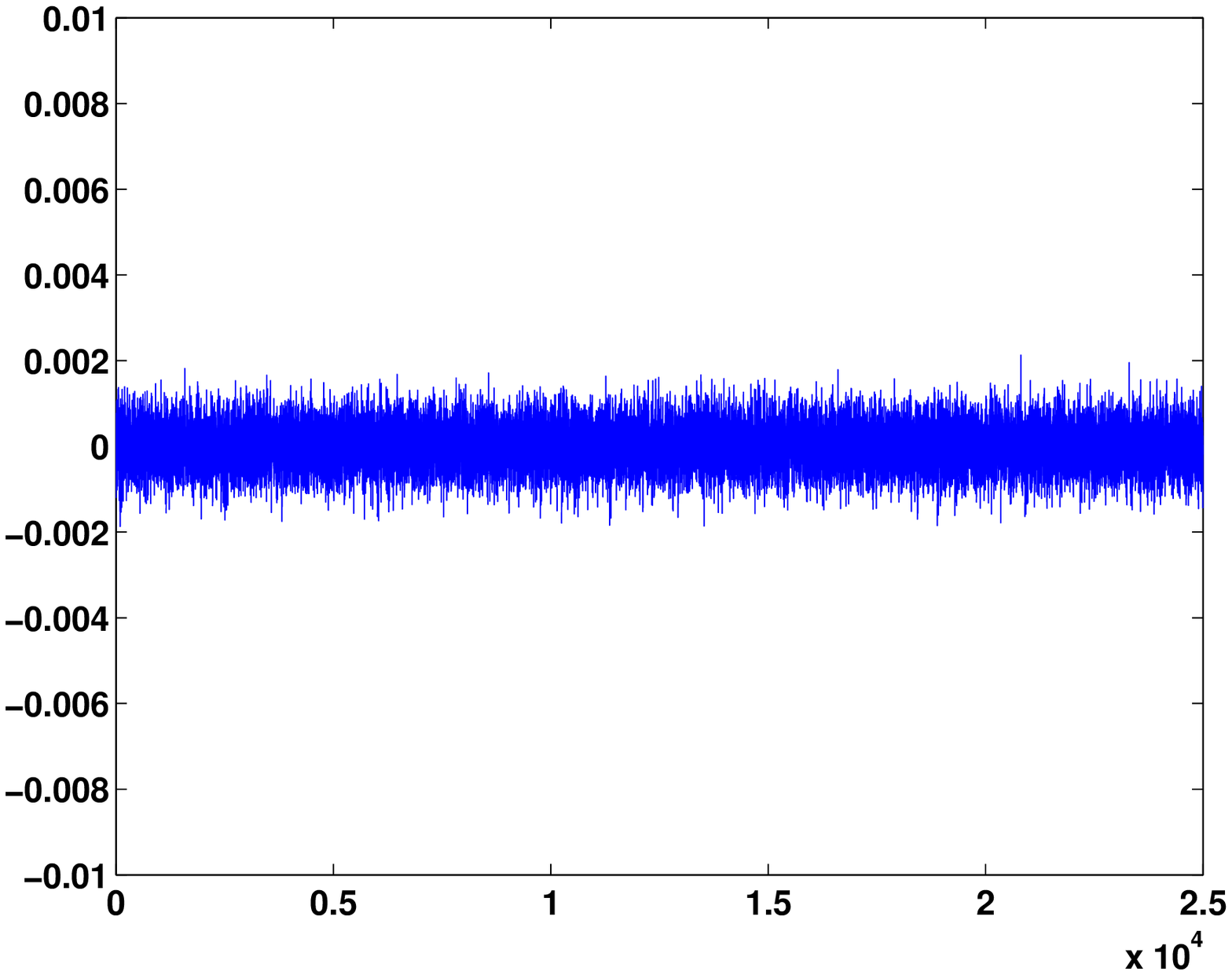} &
\includegraphics[width=2in]{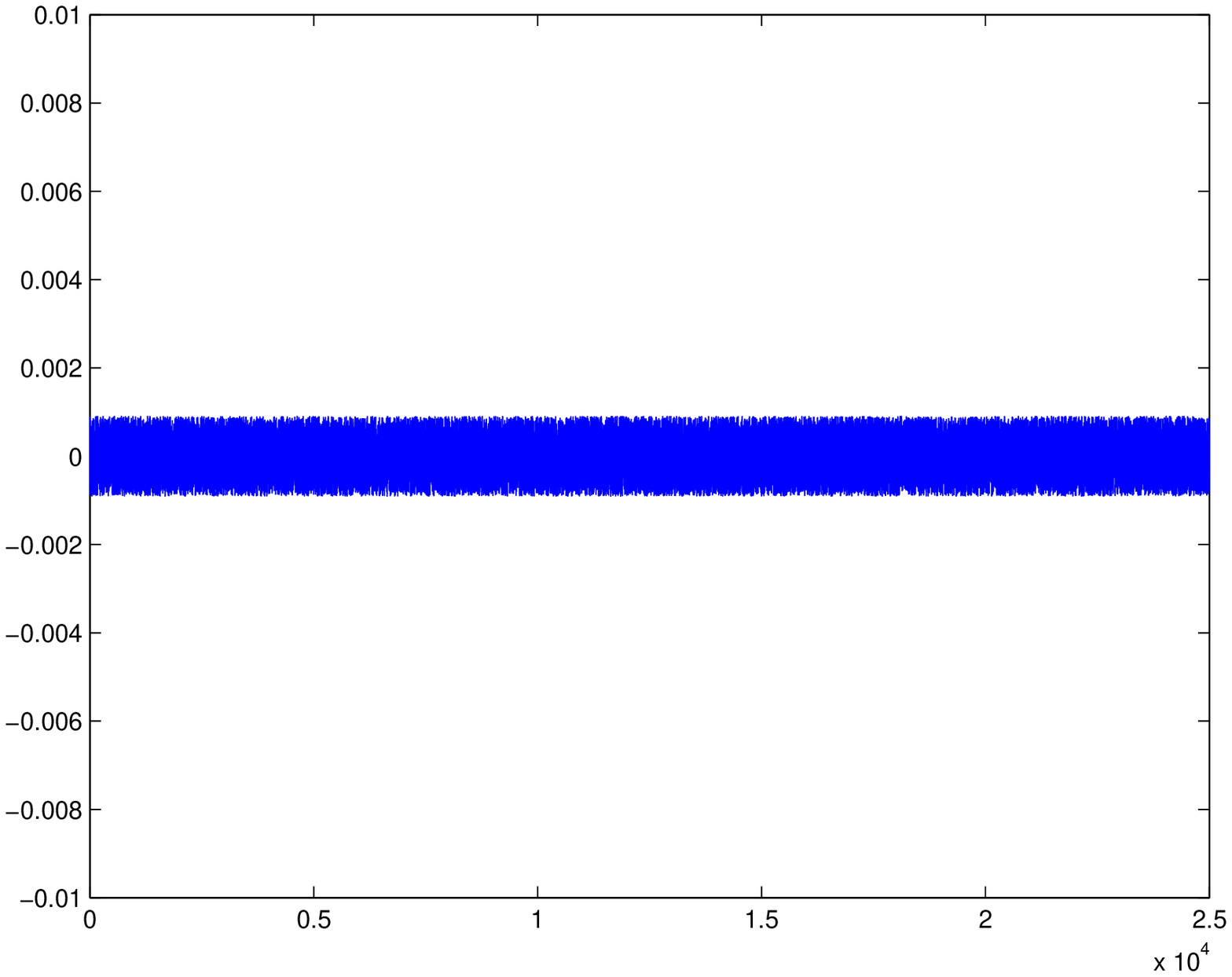} \\
(a) & (b) & (c)
\end{tabular}
}
\caption{\small\sl (a) Noiseless measurements $Ax_0$ of the {\em Boats} image.  (b) Gaussian measurement error with $\sigma = 5\cdot 10^{-4}$ in the recovery experiment summarized in the left column of 
Table~\ref{tab:image_results}.  The signal-to-noise ratio is $\|Ax_0\|_{\ell_2}/\|e\|_{\ell_2} = 4.5$.
(c) Round-off error in the recovery experiment summarized in the right column of 
Table~\ref{tab:image_results}.  The signal-to-noise ratio is $\|Ax_0\|_{\ell_2}/\|e\|_{\ell_2} = 4.3$.
}
\label{fig:measurements}
\end{figure}

We will attempt to recover the image from measurements perturbed in
two different manners.  First, as in the 1D experiments, the
measurements were corrupted by additive white Gaussian noise with
$\sigma=5\cdot10^{-4}$ so that $\sigma \cdot \sqrt{\nrow} = .0791$.
As shown in Figure~\ref{fig:measurements}, the noise level is significant; the signal-to-noise ratio is $\|Ax_0 \|_{\ell_2}/\|e\|_{\ell_2} = 4.5$.
With $\epsilon=.0798$ as in \eqref{eq:epsilon},
the recovery error is $\|\alpha^\sharp-\alpha_0\|_{\ell_2}=0.1303$ 
(the original image has unit norm).  
For comparison, the $5000$ term nonlinear approximation
error for the image is $\|\alpha_{0,5000} - \alpha_0\|_{\ell_2} =
0.050$.  Hence the recovery error is very close to the sum of the
approximation error and the size of the perturbation.

Another type of perturbation of practical interest is {\em round-off}
or {\em quantization} error.  In general, the measurements cannot be
taken with arbitrary precision, either because of limitations inherent
to the measuring device, or that we wish to communicate them using
some small number of bits.  Unlike additive white Gaussian noise,
round-off error is deterministic and signal dependent---a situation
our methodology deals with easily. 

The round-off error experiment was conducted as follows.  Using the
same scrambled Fourier ensemble, we take $25000$ measurements of {\em
  Boats}, and round (quantize) them to one digit (we restrict the
values of the measurements to be one of ten preset values, equally
spaced).  The measurement error is shown in Figure~\ref{fig:measurements}(c), and the signal-to-noise ratio is $\|Ax_0\|_{\ell_2}/\|e\|_{\ell_2} = 4.3$.
To choose $\epsilon$, we use a rough model for the size of the perturbation.
To a first approximation, the round-off error for each measurement behaves like a uniformly distributed random variable on $(-q/2,q/2)$, where 
 $q$ is the distance between quantization levels. 
Under this assumption, the size of the perturbation $\|e\|^2_{\ell_2}$ would behave like a sum of squares of uniform random variables
\[
Y = \sum_{k=1}^{\nrow} X_k^2, \qquad
X_k\sim\mathrm{Uniform}\left(-\frac{q}{2},~\frac{q}{2}\right).
\]
Here,
$\operatorname{mean}(Y) = \nrow q^2/12$ and $\operatorname{std}(Y) =
\sqrt{\nrow} q^2/[6\sqrt{5}]$.   Again, $Y$ is no larger
than $\operatorname{mean}(Y) + \lambda \operatorname{std}(Y)$ with
high probability, and we select 
\[
\epsilon^2 =  \nrow q^2/12 + \lambda \sqrt{\nrow} q^2/[6\sqrt{5}],
\]
where as before, $\lambda = 2$.  The results are summarized in the second
column of Table~\ref{tab:image_results}. 
As in the previous case, the recovery error is very close to the sum of the approximation and measurement errors.  Also note that despite the crude nature of the perturbation model, an accurate value of $\epsilon$ is chosen.

Although the errors in the recovery experiments summarized in the third row of Table~\ref{tab:image_results} are as we hoped, the recovered
images tend to contain visually displeasing high frequency oscillatory
artifacts.  To address this problem, we can solve a slightly different
optimization problem to recover the image from the same corrupted
measurements.  In place of $(P_2^\prime)$, we solve
\[
(TV) \qquad \min~\|x\|_{TV} \quad\text{subject~to}\quad \|Ax-y\|_{\ell_2}
\leq \epsilon
\]
where 
\[
\|x\|_{TV}  =  \sum_{i,j} \sqrt{(x_{i+1,j}-x_{i,j})^2 + (x_{i,j+1}-x_{i,j})^2}
  ~=~  \sum_{i,j} \left| (\nabla x)_{i,j} \right|
\]
is the {\em total variation} \cite{rudin92no} of the image $x$: the
sum of the magnitudes of the (discretized) gradient.  By substituting
$(TV)$ for $(P_2^\prime)$, we are essentially changing our model for
photograph-like images.  Instead of looking for an image with a sparse
wavelet transform that explains the observations, program $(TV)$
searches for an image with a sparse gradient (i.e. without spurious
high frequency oscillations).  In fact, it is shown in
\cite{candes04ro} that just as signals which are exactly sparse can be
recovered perfectly from a small number of measurements by solving
$(P_2)$ with $\epsilon = 0$, signals with gradients which are exactly
sparse can be recovered by solving $(TV)$ (again with $\epsilon=0$).

Figure~\ref{fig:image_example}(b) and (c) and the fourth row of
Table~\ref{tab:image_results} show the $(TV)$ recovery results.  The
reconstructions have smaller error and do not contain visually displeasing 
artifacts.

%%%%%%%%%%%%%
\begin{table}
  \caption{\small\sl Recovery results for sparse 1D signals.  
    Gaussian white noise of variance 
    $\sigma^2$ was added to each of the $\nrow=300$ measurements,
    and $(P_2)$ was solved with 
    $\epsilon$ chosen such that $\|e\|_2 \leq \epsilon$ with high probability 
    (see \eqref{eq:epsilon}).
  }
\vspace{3mm}
\centerline{
\begin{tabular}{|c||c|c|c|c|c|c|} \hline
$\sigma$ &  0.01    &    0.02   &   0.05   &   0.1   &     0.2   &    0.5 \\\hline
$\epsilon$ & 0.19 & 0.37 & 0.93 & 1.87 & 3.74 & 9.34\\\hline\hline
$\|x^\sharp - x_0\|_2$ & 0.25  &  0.49  &  1.33  &  2.55  &  4.67 &   6.61 \\\hline
\end{tabular}
}
\label{tab:sparse}
\end{table}
%%%%%%%%%%%%%

%%%%%%%%%%%%%
\begin{table}
\caption{\small\sl Recovery results for compressible 1D signals.  
Gaussian white noise of variance $\sigma^2$ was added to each measurement, and $(P_2)$ was solved with $\epsilon$ as in \eqref{eq:epsilon}.
}
\vspace{3mm}
\centerline{
\begin{tabular}{|c||c|c|c|c|c|c|} \hline
$\sigma$ &  0.01    &    0.02   &   0.05   &   0.1   &     0.2   &    0.5 \\\hline
$\epsilon$ & 0.19 & 0.37 & 0.93 & 1.87 & 3.74 & 9.34\\\hline\hline
$\|x^\sharp - x_0\|_2$ & 0.69 & 0.76 & 1.03 & 1.36 & 2.03 & 3.20 \\\hline
\end{tabular}
}
\label{tab:compress}
\end{table}
%%%%%%%%%%%%%

%%%%%%%%%%%%%
\begin{table}
  \caption{\small\sl Image recovery results.  
    Measurements of the {\em Boats} 
    image were corrupted in two different ways: 
    by adding white noise (left column) with $\sigma = 5\cdot10^{-4}$ 
    and by rounding off 
    to one digit (right column).  In each case, the image was recovered in two 
    different ways: by solving $(P_2^\prime)$ (third row) and 
    solving $(TV)$ (fourth row).  The $(TV)$ images are shown in 
    Figure~\ref{fig:image_example}.}
\vspace{3mm}
\centerline{
\begin{tabular}{|c||c|c|} \hline
 & White noise & Round-off \\\hline\hline
$\|e\|_{\ell_2}$ & 0.0789 & 0.0824 \\\hline
$\epsilon$ & 0.0798 & 0.0827 \\\hline\hline
$\|\alpha^\sharp - \alpha_0\|_{\ell_2}$ & 0.1303 & 0.1323 \\\hline
$\|\alpha^\sharp_{TV} - \alpha_0\|_{\ell_2}$ & 0.0837 & 0.0843 \\\hline
\end{tabular}
}
\label{tab:image_results}
\end{table}
%%%%%%%%%%%%%  

%%%%%%%%%%%%%
%%%%%%%%%%%%%

%%%%%%%%%%%%%
%%%%%%%%%%%%%

%%%%%%%%%%%%%
%%%%%%%%%%%%%

%----------------------------------------------------------------------
\section{Discussion}

The convex programs $(P_2)$ and $(TV)$ are simple instances of a class
of problems known as second-order cone programs (SOCP's).  As an
example, one can recast $(TV)$ as
\begin{equation}
  \label{eq:socp-tv1}
  \min \sum_{i,j} u_{i,j} \qquad \text{subject to} \qquad 
  -u_{i,j} \le \|G_{i,j} x\|_{\ell_2} \le u_{i,j}, \,\, 
\|Ax - b\|_{\ell_2} \le \epsilon, 
\end{equation}
where $G_{i,j} x = (x_{i+1,j}-x_{i,j}, x_{i,j+1}-x_{i,j})$
\cite{GoldfarbYin}. SOCP's can nowadays be solved efficiently by
interior-point methods \cite{boyd04co} and, hence, our approach is
computationally tractable.

From a certain viewpoint, recovering via $(P_2)$ is using {\em a
  priori} information about the nature of the underlying image,
i.e. that it is sparse in some known orthobasis, to overcome the
shortage of data.  In practice, we could of course use far more
sophisticated models to perform the recovery.  Obvious extensions
include looking for signals that are sparse in overcomplete wavelet or
curvelet bases, or for images that have certain geometrical structure.
The numerical experiments in Section~\ref{sec:numerical} show how
changing the model can result in a higher quality recovery {\em from
  the same set of measurements}.

%Our model assumes approximate knowledge of the noise level, i.e.
%$\|e\|_{\ell_2}$, which is the case in many practical situations. For
%example, suppose that $e$ is white noise, a vector of i.i.d.~Gaussian
%components with zero-mean and variance $\sigma^2$. Then
%$\|e\|^2_{\ell_2}/\sigma^2$ is distributed as a $\chi_\nrow^2$ with
%$\nrow$ degrees of freedom, and standard argument show that
%$\|e\|^2_{\ell_2}/\sigma^2 \le p + O(\sqrt{p})$ with overwhelming
%probability. Another possible application concerns the case where $e$
%arises because of quantization errors. Knowledge of the quantization
%step would automatically translate into appropriate estimates for
%$\|e\|_{\ell_2}$.  When the noise-level is completely unknown, one
%would have to find ways of selecting reasonable regularization
%parameters, perhaps by making some additional assumptions. This is, of
%course, is a very different subject.

%----------------------------------------------------------------------
\bibliographystyle{plain}
\bibliography{StableBib}

%----------------------------------------------------------------------
\end{document}